\documentclass[10pt,a4paper]{article}

\input epsf

\newtheorem {teor} {Theorem}

\newcommand {\ve} {\varepsilon}
\newcommand {\vp} {\varphi}

\newcommand {\av} {\bar{\cal A}}

\newcommand {\ac} {\bar{\cal C}}
\newcommand {\ad} {\bar{\cal D}}

\newcommand {\Rn} {\bf R^n}

\newcommand {\Rm} {{\bf R^m}}

\newcommand {\Rp} {\bf R^p}
\newcommand {\Rq} {\bf R^q}
\newcommand {\tve} {\left(\frac{t}{\ve}\right)}

\newcommand {\be} {\begin{equation}}
\newcommand {\en} {\end{equation}}

\newcommand {\Lm} {{\cal L}_2^m}
\newcommand {\Lq} {{\cal L}_2^q}
\newcommand {\Lp} {{\cal L}_2^p}

\begin{document}
\renewcommand{\abstractname}{}
\renewcommand{\refname}{References}
\title{Vibrational control in $H_\infty$ problems}

\author{L. A. Safonov and V. V. Strygin}

\maketitle
\begin{Large}

\begin{abstract}
We consider the application of the theory of vibrational control to 
$H_\infty$-problems. We study the possibility of introduction of
high-frequency parametric vibrations in order to decrease the minimal attainable
value of the $H_\infty$-norm. 
We prove the existence of the stabilizing solution of the Riccati equation with 
quickly oscillating coefficients. This solution is found using the averaging
technique as a series of the small parameter.
\end{abstract}

One of the goals of control in various dynamical systems is their stabilization
with respect to all or to a part of phase variables. Traditionally, the main tool 
of automatic control is feedback control, which is connected with the necessity to
measure the state of the system. However, in some systems the measurement of
the state is impossible.

As an alternative to feedback control the principle of vibrational control was
proposed by S.~M.~Meerkov \cite{M1}. The essence of vibrational control is in 
introduction of parametric vibrations into the system in order to change its dynamical 
properties. By means of vibrational control it is sometimes possible to achieve the 
desired change without measuring the system's state. As an example of such system
we can mention the inverted pendulum with the oscillating suspension point \cite{M1}. 
The stabilization of the unstable equilibrium of the pendulum in this case occurs
for sufficiently high frequency of vibrations.

In the present article we consider the application of the vibrational control theory 
[1-4] to $H_\infty$ problems [5-8].

The general $H_\infty$-problem in the linear system can be formulated as follows.
Consider the linear system
\be\label{general}
\left\{
\begin{array}{l}
\dot x=A(t)x+B_1(t)u+B_2(t)w,\\
z=Lx,
\end{array}
\right.
\en
$$  
x(0)=0,\quad t\in[0,\infty),
$$
where $x\in\Rn$, $z\in\Rm$, $u\in\Rp$, $w\in\Rq$, $A(t)$, $B_1(t)$, $B_2(t)$ 
are continuous bounded on  $[0,\infty)$ matrix functions.
We assume that the pair $(A(t),B_1(t))$ is internally stabilizable, i.e.
there exists a such continuous matrix function $K(t)$ that the system
$$
\dot x= [A(t)+B_1(t)K(t)]x
$$
is asymptotically stable.

Consider the space 
$$
{\cal L}_2^m[0,\infty)=\left\{f:[0,\infty)\to\Rm:\int_0^\infty\langle f(t),f(t)
\rangle dt<\infty\right\}
$$
and similarly defined spaces
${\cal L}_2^p[0,\infty)$, ${\cal L}_2^q[0,\infty)$.
In the sequel we shall denote these spaces as
${\cal L}_2^m$, ${\cal L}_2^p$, ${\cal L}_2^q$.

Let $U\subset\Lp$ denote the set of controls of the form $u=K(t)x$
internally stabilizing the pair $(A(t),B_1(t))$.
We shall call the quantity
$$
||\cdot||_\infty=\sup\limits_{w\in\Lq}\frac{||z||^2_{\Lm}+||u||^2_{\Lp}}{||w||^2_{\Lq}}
$$
the 
$H_\infty$-norm of the system for the given $u$.

Consider the functional
\be\label{func}
J(u,w)=||z||^2_{\Lm}+||u||^2_{\Lp}-\gamma^2||w||^2_{\Lq}.
\en

It is known that (e.g. \cite{T}) that if the Riccati equation 
\be\label{ric}
\dot R=-A^TR-RA+R\left(B_1B_1^T-\frac{1}{\gamma^2}B_2B_2^T\right)R-LL^T
\en
has the solution $R_*(t)=R_*^T(t)>0$ for $t\in[0,\infty)$ such that the system
$$
\dot x= \left[A-\left(B_1B_1^T-\frac{1}{\gamma^2}B_2B_2^T\right)R_*\right]x
$$ 
is asymptotically stable, then
\begin{enumerate}
\item The system $\dot x=[A-B_1B_1^TR_*]x$ is asymptotically stable,
\item The pair $u^*=-B_1R_*x$, $w^*=\frac{1}{\gamma^2}B_2R_*x$ is such that 
$$
J(u^*,w)=\min\limits_{u\in U}J(u,w),\ 
J(u,w^*)=\max\limits_{w\in\Lq}J(u,w),\
J(u^*,w^*)=0.
$$
\end{enumerate}

Moreover, if the solution $R_*$ with the above mentioned properties exists,
it constant if functions $A$, $B_1$, $B_2$ are constant and is periodic
if those functions are periodic. 

From 2 obviously follows that if such solution $R_*(t)$ exists, then
$J(u^*,w)\le 0$ for all $w\in\Lq$ i.e.
$$
||z||_{\Lm}\le\gamma||w||_{\Lp}
$$ 
and
$$
||\cdot||_\infty\le 0.
$$
It is also known \cite{BB} that there exists $\gamma^*(J)$ - the smallest 
of $\gamma$ for which equation (\ref{ric}) has a solution.

In the sequel we will show that $\gamma^*(J)$ can be decreased by the introduction
of vibrations into the system. Assume that the system (\ref{general}) has the form
\be\label{gen_vibr}
\left\{
\begin{array}{l}
\dot x=Ax+\frac{1}{\ve}\sin\frac{t}{\ve}Kx+B_1u+B_2w,\\
z=Lx,
\end{array}
\right.
\en
where $A$, $B_1$, $B_2$, $K$ are constant matrices and $\ve$ is a small parameter. 

As was mentioned before, the existence of a solution of the formulated above problem
is related to the existence of positive definite periodic solution of the Riccati 
equation
\be\label{ric_vibr}
\frac{dR}{dt}=-\left(A+\frac{1}{\ve}\sin\frac{t}{\ve}K\right)^TR-
R\left(A+\frac{1}{\ve}\sin\frac{t}{\ve}K\right)+
\en
$$
R\left(B_1B_1^T-\frac{1}{\gamma^2}B_2B_2^T\right)R-LL^T.
$$
Next we prove that such a solution exists for small $\ve$ if there exists
a positive definite solution of the averaged algebraic Riccati equation.

We introduce notation
$$
D=B_1B_1^T-\frac{1}{\gamma^2}B_2B_2^T,\quad C=LL^T.
$$

First, we consider the differential Riccati equation  
\be\label{ric_tau}
\frac{dR}{dt}=\ve\left[-A^T(t)R-RA(t)+RD(t)R-C(t)\right],
\en 
where $A$, $D$, $C$ are continuous $T$-periodic functions, $D$ É $C$ 
are symmetric and $\ve$ is a small parameter. 
We shall look for a solution $R_*(t)=R_*^T(t)>0$ of the equation
(\ref{ric_tau}) such that the system
$$
\dot x=\left[A(t)-D(t) R_*(t)\right]
$$
is asymptotically stable.

Consider the following spaces
\begin{itemize}
\item[] ${\bf RS^{n\times n}}$ - the space of symmetric real $n\times n$ matrices,
\item[] ${\bf CS_T^{n\times n}}$ - the space of continuous symmetric $T$-periodic
functions,
\item[] ${\bf CSZ_T^{n\times n}}$ - subspace of 
${\bf CS_T^{n\times n}}$, consisting of functions with the zero average.
\end{itemize}

Obviously, 
$$
{\bf CS_T^{n\times n}}={\bf RS^{n\times n}}\oplus 
{\bf CSZ_T^{n\times n}},
$$
if the elements of ${\bf RS^{n\times n}}$ are considered as constant functions.

Now we introduce the averaging operator 
$M_0:{\bf CS_T^{n\times n}}\to{\bf RS^{n\times n}}$ defined by the equality
$$
M_0F=\frac{1}{T}\int_0^TF(t)dt.
$$
Clearly, $M_0$ is a projector to the subspace ${\bf RS^{n\times n}}$. 
We also define the operator
$T_0$ on ${\bf CSZ_T^{n\times n}}$ 
$$
T_0F=F-M_0F.
$$
For brevity we will denote $M_0F$ as $\bar F$.

Consider the algebraic Riccati equation
\be\label{ric_aver}
-\bar A^TR-R\bar A+R\bar DR-\bar C=0.
\en

\begin{teor}
If equation (\ref{ric_aver}) has such solution $R_0=R_0^T>0$ that 
the matrix $\bar A-\bar DR_0$ is Hurwitz, than there exists such $\ve_0>0$ that
for any $\ve<\ve_0$ equation (\ref{ric_tau}) has such a $T$-periodic
symmetric positive definite solution $R_*(t)$, that the system 
$$
\dot x=[A(t)-D(t)R_*(t)]x
$$
is asymptotically stable.

This solution can be approximately found as
$$
R_{(N)}(t)=R_0+\ve(R_1+\Pi_1(t))+\ldots+
\ve^N(R_N+\Pi_N(t))+\ve^{N+1}\Pi_{N+1}(t),
$$  
where $R_0,\ldots,R_N\in{\bf RS^{n\times n}}$, 
$\Pi_1(t),\ldots,\Pi_{N+1}(t)\in{\bf CSZ_T^{n\times n}}$,
and for some $C>0$ independent of $\ve$
\be\label{Rest}
\sup\limits_{t\in[0,T]}||R_*(t)-R_{(N)}(t)||<C\ve^{N+1}.
\en
\end{teor}
 
{\it Proof.} Formally substituting $R_{(N)}$ into Eq. 
(\ref{ric_tau}), we have
$$
\ve\dot\Pi_1+\ve^2\dot\Pi_2+\ldots+\ve^{N+1}\dot\Pi_{N+1}=
$$
$$
\ve\left[-A^T(t)\left(R_0+\ve(R_1+\Pi_1(t))+\ldots\right)-\right.
$$
$$
\left(R_0+\ve(R_1+\Pi_1(t))+\ldots\right)A(t)+
$$
\be\label{subs}
\left.\left(R_0+\ve(R_1+\Pi_1(t))+\ldots\right)D(t)
\left(R_0+\ve(R_1+\Pi_1(t))+\ldots\right)-C(t)
\right].
\en
Equating in (\ref{subs}) the terms of order $\ve^1$ we have
$$
\frac{d\Pi_1}{dt}=-A^T(t)R_0-R_0A(t)+R_0D(t)R_0-C(t).
$$
Averaging the latter equation yields
\be\label{ric_aver1}
-\bar A^TR_0-R_0\bar A+R_0\bar DR_0-\bar C=0
\en
É
$$
\Pi_1(t)=\int_0^t T_0\left[-A^T(s)R_0-R_0A(s)+R_0D(s)R_0-C(s)\right]ds.
$$
According to the assumptions of the theorem, equation (\ref{ric_aver1}) has a
positive definite solution $R_0\in{\bf RS^{n\times n}}$, for which the matrix
$\bar A-\bar DR_0$ is Hurwitz.

Equating the terms of order $\ve^2$, we have
$$
\frac{d\Pi_2}{dt}=-A^T(t)(R_1+\Pi_1(t))-(R_1+\Pi_1(t))A(t)+
$$
$$
(R_1+\Pi_1(t))D(t)R_0+R_0D(t)(R_1+\Pi_1(t)).
$$
Averaging of the latter equality yields
\be\label{ric_aver2}
(\bar A^T-\bar DR_0)^TR_1+R_1(\bar A^T-\bar DR_0)+
\en
$$
M_0[-A^T(t)\Pi_1(t)-\Pi_1(t)A(t)+\Pi_1(t)D(t)R_0+
R_0D(t)\Pi_1(t)]=0.
$$
Because the matrix $\bar A-\bar DR_0$ is Hurwitz, the Lyapunov
equation (\ref{ric_aver2}) has the unique solution $R_1$.

The matrix $\Pi_2(t)$ is defined by 
$$
\Pi_2(t)=\int_0^t T_0\left[A^T(s)(R_1+\Pi_1(s))-(R_1+\Pi_1(s))A(s)+
\right.
$$
$$
\left.
(R_1+\Pi_1(s))D(s)R_0+R_0D(s)(R_1+\Pi_1(s))\right]ds.
$$

Similarly we can find $R_2,\ldots,R_N$ and $\Pi_3(t),\ldots,\Pi_{N+1}(t)$.

Next we prove the existence of solution $R_*(t)$, satisfying the
conditions of the present theorem. Obviously,
$$
\frac{dR_{(N)}}{dt}=\ve\left[-A^T(t)R_{(N)}-R_{(N)}A(t)+
R_{(N)}D(t)R_{(N)}-C(t)+\right.
$$
$$
\left.\ve^{N+1}g(t,\ve)\right],
$$
where $g(t,\ve)$ - $T$ is a periodic in $t$ and uniformly bounded in
$\ve\in(0,\ve_0]$ function.

First we prove that the equation
\be\label{dtheta}
\frac{d\Theta}{dt}=\ve\left[-(A-DR_{(N)})^T\Theta-\Theta(A-DR_{(N)})+
\Theta D\Theta-\ve^{N+1}g(t,\ve)\right]
\en
has a $T$-periodic solution.
Note that the function $H(\Theta)=\Theta D\Theta$ is such that for any 
$\sigma>0$ there exists $\delta>0$ such that if $||\Theta_1||<\delta$ and
$||\Theta_1||<\delta$, than
$$
||H(\Theta_1-\Theta_2)||<\sigma||\Theta_1-\Theta_2||.
$$

Let $\Gamma=M_0\Theta$ and $\Delta=T_0\Theta$. Than applying operators
$M_0$ and $T_0$ to the left and right hand sides of Eq.
(\ref{dtheta}), we have
\be\label{dgamma}
0=\int_0^T\left[-(A-DR_{(N)})^T(\Gamma+\Delta)-(\Gamma+\Delta)(A-DR_{(N)})+
\right.
\en
$$
\left.
(A-DR_{(N)})D(A-DR_{(N)})-\ve^{N+1}g(t,\ve)\right]dt,
$$
\be\label{ddelta}
\dot\Delta=\ve T_0\left[-(A-DR_{(N)})^T(\Gamma+\Delta)-
(\Gamma+\Delta)(A-DR_{(N)})+\right.
\en
$$
\left.
(A-DR_{(N)})D(A-DR_{(N)})-
\ve^{N+1}g(t,\ve)\right].
$$
From Eq. (\ref{dgamma}) we have
\be\label{AB}
{\bf A}\Gamma+{\bf B}\Delta+G_1(\Gamma,\Delta,\ve)=0,
\en
where the operators ${\bf A}:{\bf RS^{n\times n}}\to{\bf RS^{n\times n}}$ and
${\bf B}:{\bf CSZ_T^{n\times n}}\to{\bf RS^{n\times n}}$
are defined by the equality
$$
{\bf A}\Gamma=-(\bar A-\bar DR_0)^T\Gamma-\Gamma(\bar A-\bar DR_0),
$$
$$
{\bf B}\Delta=-M_0\left[(A(t)-D(t)R_0)^T\Delta+
\Delta(A(t)-D(t)R_0)\right],
$$
and
$$
G_1(\Gamma,\Delta,\ve)=M_0\left[\ve D\rho_1(t,\ve)(\Gamma+\Delta)+\right.
$$
$$
\left.
\ve(\Gamma+\Delta)D\rho_1(t,\ve)+(\Gamma+\Delta)D(\Gamma+\Delta)-
\ve^{N+1}g(t,\ve)\right],
$$
where $\rho_1(t,\ve)=(R_{(N)}(t)-R_0)/\ve$.

Obviously, for any $\sigma_1>0$ there exist $\ve_0$ and $\delta$ such that 
if $||\Gamma_1||_{{\bf RS}}<\delta$, $||\Gamma_2||_{{\bf RS}}<\delta$, 
$||\Delta_1||_{{\bf CSZ}}<\delta$, $||\Delta_2||_{{\bf CSZ}}<\delta$ and $\ve<\ve_0$ than
\be\label{lipsh1}
||G_1(\Gamma_1,\Delta_1,\ve)-G_1(\Gamma_2,\Delta_2,\ve)||_{{\bf CSZ}}<
\sigma(||\Gamma_1-\Gamma_2||_{{\bf RS}}+||\Delta_1-\Delta_2||_{{\bf CSZ}}).
\en
Since the matrix $\bar A-\bar DR_0$ is Hurwitz, the operator ${\bf A}$ is
invertible and from (\ref{AB}) follows that
\be\label{gammaeq}
\Gamma=-{\bf A}^{-1}{\bf B}\Delta-{\bf A}^{-1}G_1(\Gamma,\Delta,\ve).
\en
From Eq. (\ref{ddelta}) we have 
\be\label{deltaeq}
\Delta=\ve G_2(\Gamma,\Delta,\ve),
\en
where
$$
G_2(\Gamma,\Delta,\ve)=
$$
$$
\int_0^t T_0\left[-(A(s)-D(s)R_{(N)}(s))^T(\Gamma+\Delta(s))-\right.
$$
$$
(\Gamma+\Delta(s))(A(s)-D(s)R_{(N)}(s))+
$$
$$
\left.
(\Gamma+\Delta(s))D(\Gamma+\Delta(s))-\ve^{N+1}g(s,\ve)\right]ds.
$$
$G_2(\Gamma,\Delta,\ve)$ satisfy a condition similar to (\ref{lipsh1}), 
namely that for any $\sigma>0$ there exist such $\ve_0$ and $\delta$ that 
if $||\Gamma_1||_{{\bf RS}}<\delta$, $||\Gamma_2||_{{\bf RS}}<\delta$, 
$||\Delta_1||_{{\bf CSZ}}<\delta$,
$||\Delta_2||_{{\bf CSZ}}<\delta$ and $\ve<\ve_0$ than
\be\label{lipsh2}
||G_2(\Gamma_1,\Delta_1,\ve)-G_2(\Gamma_2,\Delta_2,\ve)||_{{\bf CSZ}}<
\sigma(||\Gamma_1-\Gamma_2||_{{\bf RS}}+||\Delta_1-\Delta_2||_{{\bf CSZ}}).
\en     

System (\ref{gammaeq}), (\ref{deltaeq}) can be rewritten as
\be\label{gd}
\left(
\begin{array}{c}
\Gamma\\
\Delta
\end{array}
\right)
=
\left(
\begin{array}{cc}
0&-{\bf A}^{-1}{\bf B}\\
0&0
\end{array}
\right)
\left(
\begin{array}{c}
\Gamma\\
\Delta
\end{array}
\right)
+
\left(
\begin{array}{c}
-{\bf A}^{-1}G_1(\Gamma,\Delta,\ve)\\
\ve G_2(\Gamma,\Delta,\ve)
\end{array}
\right)
\en
or
$$
(\Gamma,\Delta)=\vp(\Gamma,\Delta),
$$ 
where
$$
\vp(\Gamma,\Delta)=
\left(
\begin{array}{cc}
0&-{\bf A}^{-1}{\bf B}\\
0&0
\end{array}
\right)
\left(
\begin{array}{c}
\Gamma\\
\Delta
\end{array}
\right)
+
\left(
\begin{array}{c}
-{\bf A}^{-1}G_1(\Gamma,\Delta,\ve)\\
\ve G_2(\Gamma,\Delta,\ve)
\end{array}
\right).
$$
The map $\vp$ acts in the space
$$
{\bf X}={\bf RS^{n\times n}}\times{\bf CSZ_T^{n\times n}},
$$
where we introduce the vector norm
$$
||(\Gamma,\Delta)||_{\bf X}=(||\Gamma||_{\bf RS},||\Delta||_{\bf CSZ})^T.
$$
Next we show that there exists a ball in ${\bf X}$, within which the map
$\vp$ satisfies the generalized contraction principle \cite{Kras}, i.e.
there exists the matrix $K\in{\bf R^{2\times 2}}$ with non-negative elements and 
the spectral radius $\rho(K)<1$ such that
$$
||\vp(\Gamma_1,\Delta_1)-\vp(\Gamma_2,\Delta_2)||_{\bf X}
\le K||(\Gamma_1,\Delta_1)-(\Gamma_2,\Delta_2)||_{\bf X},
$$
where the relation of inequality is understood in the following sense:

it is said that 
$x\le y$ ($x,y\in{\bf R^2}$) if $x_i\le y_i$ for $i=1,2$.

From (\ref{gd}) and estimates (\ref{lipsh1}), (\ref{lipsh2}) follows that if  
$||\Gamma_1||_{{\bf RS}}<\delta$, $||\Gamma_2||_{{\bf RS}}<\delta$, 
$||\Delta_1||_{{\bf CSZ}}<\delta$, $||\Delta_2||_{{\bf CSZ}}<\delta$,
then
\be\label{contr}
||\vp(\Gamma_1,\Delta_1)-\vp(\Gamma_2,\Delta_2)||_{\bf X}\le
\en
$$
\left(
\begin{array}{cc}
||{\bf A}^{-1}||\sigma_1&||{\bf A}^{-1}{\bf B}||+
||{\bf A}^{-1}||\sigma_1\\
\ve\sigma_2&\ve\sigma_2
\end{array}
\right)
||(\Gamma_1,\Delta_1)-(\Gamma_2,\Delta_2)||_{\bf X},
$$
where $\sigma_1,\sigma_2\to 0$ for $\delta\to 0$.

Therefore, the map $\vp$ satisfies the generalized contraction principle for
sufficiently small $\delta$. Let us now find the value of $\delta$ for which
$\vp$ maps the ball of the radius $\delta$ in the space ${\bf X}$ into itself. 

From (\ref{lipsh1}) and (\ref{lipsh2}) follows that
$$
||G_1(\Gamma,\Delta,\ve)-G_1(0,0,\ve)||\le\sigma_1(||\Gamma||+||\Delta||),
$$
$$
||G_2(\Gamma,\Delta,\ve)-G_2(0,0,\ve)||\le\sigma_2(||\Gamma||+||\Delta||),
$$
for $||\Gamma||<\delta$ É $||\Delta||<\delta$. This implies that 
$$
||G_1(\Gamma,\Delta,\ve)||<C_1\ve^{N+1}+\sigma_1(||\Gamma||+||\Delta||),
$$
$$
||G_2(\Gamma,\Delta,\ve)||<C_2\ve^{N+1}+\sigma_2(||\Gamma||+||\Delta||),
$$
where $C_1$, $C_2$ are some constants.

From (\ref{contr}) follows that
$$
||\vp(\Gamma,\Delta)||_{\bf X}\le
$$
$$
\left(
\begin{array}{cc}
||{\bf A}^{-1}||\sigma_1&||{\bf A}^{-1}{\bf B}||+||{\bf A}^{-1}||\sigma_1\\
\ve\sigma_2&\ve\sigma_2
\end{array}
\right)
||(\Gamma,\Delta)||_{\bf X}+
\left(
\begin{array}{c}
C_1\ve^{N+1}\\
C_2\ve^{N+1}
\end{array}
\right),
$$
i.e. for some $C^*$ $\vp$ maps the ball of the radius $\rho=C^*\ve^{N+1}$ in the
space ${\bf X}$ into itself and is contracting, and hence has in this ball the
unique fixed point $(\Gamma_*,\Delta_*(t))$ such that
$\Theta_*(t)=\Gamma_*+\Delta_*(t)$ is the solution of Eq.
(\ref{dtheta}). It is easy to verify that $R_*(t)=R_{(N)}(t)+\Theta(t)$ 
is a solution of equation(\ref{ric_tau}).  

Next we prove the asymptotic stability of the system 
\be\label{asymp}
\dot x=[A(t)-D(t)R_*(t)]x.
\en
Substituting $R_*(t)$ in the form $R_0+\ve\rho_1(t)+\Theta(t)$
into (\ref{asymp}) and taking into account that the matrix $\bar A-\bar DR_0$
is Hurwitz, we have that (\ref{asymp}) is asymptotically stable for small $\ve$.

Positive definiteness of $R_*(t)$ for small $\ve$ follows from the positive 
definiteness of $R_0$.

Theorem 1 is proved.

Consider Eq. (\ref{ric_vibr}). Introduce the fast time $\tau=t/\ve$. Then
\be\label{ric_tau1}
\frac{dR}{d\tau}=-\sin\tau(K^TR+RK)+\ve(-A^TR-RA+RDR-C).
\en

Putting $\ve=0$, we have
\be\label{ric_tau_lin}
\frac{dR}{d\tau}=-\sin\tau(K^TR+RK).
\en
Equation (\ref{ric_tau_lin}) has the general solution 
\be\label{P}
R=\Psi(\tau)P\Psi^T(\tau),
\en
where $P$ is an arbitrary constant, $\Psi(\tau)=\exp(K^T(\cos\tau-1))$.

Introducing $P$, defined by (\ref{P}), as the new variable in (\ref{ric_tau1}),
we have
\be\label{ric_eps}
\frac{dP}{d\tau}=\ve\left(-{\cal A}^T(\tau)P-P{\cal A}(\tau)+
P{\cal D}(\tau)P-{\cal C}\right),
\en
where 
$$
{\cal A}(\tau)=\Psi^T(\tau)A(\Psi^T)^{-1}(\tau)\quad
{\cal D}(\tau)=\Psi^T(\tau)D\Psi(\tau),
$$
$$
{\cal C}(\tau)=\Psi^{-1}(\tau)C(\Psi^T)^{-1}(\tau).
$$
Consider the algebraic Riccati equation
\be\label{ric_aver3}
-\av^TP-P\av+P\ad P-\ac=0,
\en
which is the averaging of the right hand side of (\ref{ric_eps}).
\begin{teor}
If equation (\ref{ric_aver3}) has a solution $P_0=P_0^T>0$ such that
the matrix $\av-\ad P_0$ is Hurwitz, then there exists a $\ve_0>0$ such
that for any $\ve<\ve_0$ equation (\ref{ric_vibr}) has a $2\pi\ve$-periodic
symmetric positive definite solution $R_*(\frac{t}{\ve})$ for which the zero
solution of the system
$$
\frac{dx}{dt}=\left[{\cal A}\tve-
{\cal D}\tve R_*\tve\right]x
$$
is exponentially stable.

For any $N\ge 0$ this solution can be found in the form
$$
R_{(N)}\tve=
\Psi\left(\frac{t}{\ve}\right)
\left[P_0+\ve\left(P_1+\Pi_1\tve\right)+
\ldots+\right.
$$
$$
\left.
\ve^N\left(P_N+\Pi_N\tve\right)+
\ve^{N+1}\Pi_{N+1}\tve\right] 
\Psi^T\left(\frac{t}{\ve}\right),
$$
where $P_1,\ldots,P_N\in{\bf RS^{n\times n}}$, 
$\Pi_1,\ldots,\Pi_{N+1}\in{\bf CSZ_{T\ve}^{n\times n}}$,
and for some $C>0$
$$
\sup\limits_{\tau\in[0,T\ve]}
\bigg|\bigg|R_*\tve-R_{(N)}\tve\bigg|\bigg|<C\ve^{N+1}.
$$
\end{teor} 

{\it Proof.}
According to Theorem 1, equation (\ref{ric_eps}) has a positive definite 
stabilizing solution $P_*(\tau)$ such that for 
$$
P_{(N)}(\tau)=P_0+\ve(P_1+\Pi_1(\tau))+\ldots+\ve^N(P_N+\Pi_N(\tau))+
\ve^{N+1}\Pi_{N+1}(\tau)
$$  
the estimate (\ref{Rest}) is satisfied. It is obvious that
$R(t/\ve)=\Psi(t/\ve) P_*(t/\ve)\Psi^T(t/\ve)$ is a solution of (\ref{ric_vibr}) and
for $R_{(N)}(t/\ve)$ and $R_*(t/\ve)$ a similar estimate is satisfied, and that the 
system 
$$
\dot x=[{\cal A}(\tau)-{\cal D}(\tau)R_*(\tau)]x
$$
is exponentially stable.

{\bf Example.} Let matrices $A$, $B$, $C$ and $K$ have the form
$$
A=
\left(
\begin{array}{cc}
0&1\\
-0.27&-2.8
\end{array}
\right),\quad
B=
\left(
\begin{array}{c}
0\\
1
\end{array}
\right),
$$
$$
C=
\left(
\begin{array}{cc}
1&0\\
0&1
\end{array}
\right),\quad
K=
\left(
\begin{array}{cc}
0&0\\
0&k
\end{array}
\right).
$$
The averaged Riccati equation (\ref{ric_aver3}) in this case has the form 
$$
\left(
\begin{array}{cc}
0&-0,27-k^2/2\\
1&-2,8
\end{array}
\right)P
+
P
\left(
\begin{array}{cc}
0&1\\
-0,27-k^2/2&-2,8
\end{array}
\right)+
$$
$$
{1\over{\gamma^2}}
P
\left(
\begin{array}{cc}
0&0\\
0&1
\end{array}
\right)
P+
\left(
\begin{array}{cc}
1&0\\
0&1
\end{array}
\right)=0.
$$
The table below shows the values $\gamma^*_K$ for different $k$.
These values are found numerically using the conditions for the existence
of solutions of Riccati equations, given in \cite{Hewer}.

\bigskip
\begin{tabular}{|p{2cm}|p{2cm}|}
\hline
$k$&$\gamma^*_K$\\
\hline
0&3,704\\
0,25&3,320\\
0,5&2,532\\
0,75&1,815\\
1&1,300\\
1,25&0,925\\
1,5&0,717\\
1,75&0,556\\
\hline
\end{tabular}
\bigskip

This example shows the possibility of introduction of vibrations into a
linear system, which decrease $\gamma^*_K$ i.e. decrease the influence
of disturbances on the output of the system.

\end{Large}

\end{document}